\documentclass[graybox]{svmult}

\usepackage{type1cm}        % activate if the above 3 fonts are
\usepackage{makeidx}         % allows index generation
\usepackage{graphicx}        % standard LaTeX graphics tool
\usepackage{multicol}        % used for the two-column index
\usepackage[bottom]{footmisc}% places footnotes at page bottom

\usepackage{newtxtext}       % 
\usepackage[varvw]{newtxmath}       % selects Times Roman as basic font

\makeindex             % used for the subject index
\def\og{\leavevmode\raise.3ex\hbox{$\scriptscriptstyle\langle\!\langle\,$}}

\def\fgf{\/\leavevmode\raise.3ex\hbox{$\scriptscriptstyle\,\rangle\!\rangle$}}

\def\fg{\fgf\ }

\def\Demd#1|{\parindent=0pt\par{\sl D\'emonstration d#1}\pointir\parindent=20pt}

\let\petcap\sc

\def\finc{\vskip12pt}

\def\Dem{\parindent=0pt\par{\sl D\'emonstration}\pointir\parindent=20pt}

\def\Demo#1|{\parindent=0pt\par{\sl D\'emonstration #1}\pointir\parindent=20pt}

\makeatletter
\newdimen\indentTh\indentTh=0pt
\def\p@int{{\rm .}}
\def\p@intir{\discretionary{\rm .}{}{\rm .\kern.35em---\kern.7em}}
\def\pointir{\afterassignment\pointir@\let\next=}
\def\pointir@{\ifx\next\par\p@int\else\p@intir\fi\next}
\long\def\Thc#1|#2\finc{\Th{}{#1}{\pointir}{#2}}
\long\def\Th#1#2#3#4{\parindent=\indentTh\par\vskip5pt
{#1}{\petcap #2}{\sl #3}\parindent=20pt{\sl #4\par}\vskip 5pt\parindent=20pt}

\newdimen\indentssec\indentssec=20pt
\newdimen\indentrem\indentssec=0pt

\long\def\Remarque#1#2#3#4{\parindent=\indentrem\par\vskip5pt
{#1}{ \sl #2}{\sl #3}\parindent=20pt{#4\par}
\vskip 5pt\parindent=20pt}
\long\def\Rmc#1|#2\finc{\Remarque{}{#1}{\pointir}{#2}}

\long\def\rmc#1|#2\finc{\remarque{}{#1}{\pointir}{#2}}
\long\def\parc#1\finc{\remarque{}{}{}{#1}}
\long\def\remarque#1#2#3#4{\parindent=\indentrem\par\vskip5pt
{#1}{\sl #2}{\sl #3}\parindent=20pt{#4\par}
\vskip 5pt\parindent=20pt}
\makeatother

\def\Demdsp#1|{\parindent=0pt\par{\sl D\'emonstration d#1.}\parindent=20pt}

\def\hdfl#1#2{\smash{\mathop{\hbox to 12mm{\rightarrowfill}}
\limits^{\scriptstyle#1}_{\scriptstyle#2}}}

\def\hdhfl#1#2{\smash{\mathop{\hbox to 12mm{\hookrightarrowfill}}
\limits^{\scriptstyle#1}_{\scriptstyle#2}}}

\def\hgfl#1#2{\smash{\mathop{\hbox to 12mm{\leftarrowfill}}
\limits^{\scriptstyle#1}_{\scriptstyle#2}}}

\def\hghfl#1#2{\smash{\mathop{\hbox to 12mm{\hookleftarrowfill}}
\limits^{\scriptstyle#1}_{\scriptstyle#2}}}

\begin{document}

%%%%% d0.tex %%%%%%
 
\title*{Une remarque sur l'invariant de Arf%
.}

\author{{\it par} Alexis Marin}
\maketitle

\centerline{\bf R\'esum\'e}

Utilisant le \og formalisme de Witt\fg  d\'ecrit dans {\bf [2]} on donne
une construction, sans calcul de changement de base, de l'invariant de Arf
des forme quadratiques en caract\'eristique 2.
\vskip 5mm

\centerline{\bf Abstract}

An alternative construction of the Arf-invariant of quadratic forms
in characteristic 2.
\finc
\vskip12pt

\noindent 2020 {\sl Mathematics Subject Classification}~:\hfill\break
Primary: 11E81; Secondary: 11E04.

\def\keywordname{{\bf Key Words and Phrases:}}%

{\parindent=0pt\par\vskip .3cm
\vskip 0mm plus -20mm minus 1,5mm\penalty-50
{\bf 1.\ \/}%
{\petcap Conoyau d'Artin-Schreier
et extensions radicielles}{\bf \pointir}%
\nobreak\parindent=20pt}%
Soient $k$ un corps de caract\'eristique  $p$ positive,
$\Phi_k : k\rightarrow k, x\mapsto x^p$ et ${\cal P\/}_k$ ses Frobenius  
et {\it Artin-Schreier\/}
({\it i.e.\/} l'endomorphisme ${\cal P\/}=\Phi-{\rm Id\/} : k\rightarrow k$
du groupe additif de $k$).

\Thc Lemme 0| Une extension radicielle $\iota : k\subset K$ 
induit un isomorphisme
$$[\iota] : k/{\cal P}_k\,(k)
\xrightarrow[]{\sim} K/{\cal P}_K\,(K)$$
\finc
\Dem Soit $a$ dans $k$ avec $a={\cal P\/}_K\,(x)=x^p-x$
o\`u $x$ est dans $K$. Alors $x$ est dans $k$,
puisque l'\'equation $X^p-X-a=0$
est s\'eparable, d'o\`u l'injectivit\'e.

Soit $x$ dans $K$ de hauteur $n$ sur $k$, ainsi
$x^{p^{n}}=y$ est dans $k$, et
$x=y-x^{p^n}+x^{p^{n-1}}+\cdots -x^{p}+x
= y-{\cal P\/}_K\,(x^{p^{n-1}}+\cdots+x)$,
d'o\`u la surjectivit\'e.
{\parindent=0pt\par\vskip .3cm
\vskip 0mm plus -20mm minus 1,5mm\penalty-50
{\bf 2.\ \/}%
{\petcap Vocabulaire du formalisme de Witt\/} (Cf. {\bf [2]}){\bf \pointir}%
\nobreak\parindent=20pt}%
Une {\it forme quadratique\/} sur un corps $k$ est une application
$q : V\rightarrow k$ d\'efinie sur un $k\/$-espace vectoriel $V$
de dimension finie
et telle que $(x,\, y)\mapsto b\,(x,\,y)=q\,(x+y)-q\,(x)-q\,(y)$
est bilin\'eaire non d\'eg\'en\'er\'ee, la
{\it forme bilin\'eaire\/} de $q$.
Un $b\/$-{\it Lagrangien\/} de $b$ est un sous-espace $L$ de $V$
\'egal \`a son orthogonal. C'est un $q\/$-{\it Lagrangien\/}
si de plus $q\,(L)=0$. Une forme quadratique est {\it neutre\/} si
elle a un $q\/$-Lagrangien.
Le quotient du mono\"{\i}de de somme
orthogonale des formes quadratiques sur $k$
par le sous-mono\"{\i}de des formes neutres est un groupe, le
{\it groupe de Witt quadratique\/} $WQ\,(k)$ du corps $k$.

{\parindent=0pt\par\vskip .3cm
\vskip 0mm plus -20mm minus 1,5mm\penalty-50
{\bf 3.\ \/}%
{\petcap Formes quadratiques sur un corps parfait
de caract\'eristique $2$\/}{\bf \pointir}%
\nobreak\parindent=20pt}%
Soit $q : V\!\rightarrow\!K$ une forme quadratique
sur un corps parfait de caract\'eristique $2$. On note $\sqrt{^{\ }}$ l'inverse
de $\Phi_K$. La
forme bilin\'eaire $b$ de $q$ est altern\'ee car
$b\,(x, x)=2\,q\,(x)=0$
pour tout $x$ de $V$. Ayant une base symplectique (C.f. {\bf [3]} I 3.5),
cette forme  a un $b\/$-{\it Lagrangien\/}, $L$. La fonction
$\sqrt{q}$
est lin\'eaire sur $L$
car pour $x, y$ de $L=L^\perp$ et $\lambda$ de $K$  on a :
$\sqrt{q\,(x+y)}=\sqrt{q\,(x)+b\,(x, y)+q\,(y)}
=\sqrt{q\,(x)+q\,(y)}
=\sqrt{q\,(x)}+\sqrt{q\,(y)}$
 et $\sqrt{q\,(\lambda\, x)}=
\sqrt{\lambda^2\, q(x)}
=\lambda\,\sqrt{q\,(x)}$.
Il y a donc dans $V$, car $b$ est non d\'eg\'en\'er\'ee,
un {\it vecteur de Wu\/} $\omega$ {\it de\/} $L$, 
bien d\'efini  modulo $L^\perp=L$ par la relation
$b\,(\omega, l)=-\sqrt{q\,(l)}$ pour tout
$l$ de $L$. Donc
$q\,(\omega+l)=q\,(\omega)+b\,(\omega, l)+q\,(l)=
q\,(\omega)-\sqrt{q\,(l)}+\Phi\,(\sqrt{q\,(l)})=
q\,(\omega)+{\cal P\/}\,(\sqrt{q\,(l)})$ :
\Thc Lemme 1|Quand $\omega+l$ parcourre la classe de $\omega$
modulo $L$ alors $q\,(\omega+l)$ est nul si $q$ est nulle sur $L$
 et d\'ecrit une classe modulo ${\cal P\/}\,(K)$
sinon.
\finc
\Thc Proposition| La classe modulo ${\cal P\/}\,(K)$ de
$q\,(\omega)$ ne d\'epend que de la classe de Witt de $q$ et
$[q]\mapsto [q\,(\omega)]$ induit l'isomorphisme
$${\rm Parf\/}: WQ\,(K)
\xrightarrow[]{\sim} K/{\cal P\/}\,(K)$$
du groupe de Witt quadratique du corps $K$ sur
le conoyau de son Artin-Schreier.
\finc
\Dem Si $\omega$ est vecteur de Wu d'un $q\/$-Lagrangien
alors $q\,(\omega)=0$. De plus, si $\omega_i$, pour $i=1, 2$,
sont vecteurs de Wu  
de $b\/$-Lagrangiens $L_i$ de formes quadratiques $q_i$,
alors $\omega=\omega_1+\omega_2$ est vecteur de Wu du $b\/$-Lagrangien
 $L_1\oplus L_2$ de la somme orthogonale $q=q_1\oplus q_2$
avec $q\,(\omega)=q_1\,(\omega)+q_2\,(\omega)$.
Ainsi pour que ${\rm Parf\/}$
soit un morphisme bien d\'efini sur le groupe de Witt quadratique il suffit,
d'apr\`es le {\petcap Lemme 1\/},
de montrer que deux $b\/$-Lagrangiens $L_1$ et $L_2$ d'une
forme quadratique $q$ ont un vecteur de
Wu commun : Soit $\omega_i\/$ vecteur de Wu de $L_i$.
Pour $l$ dans $L_1\cap L_2$ on a
$b\,(\omega_1-\omega_2,\,l)=q\,(l)-q\,(l)=0$. Ainsi $\omega_1-\omega_2$,
\'etant
dans  l'orthogonal $(L_1\cap L_2)^\perp=L_1^\perp+L_2^\perp$, s'\'ecrit
$\omega_1-\omega_2=m_1+m_2$ o\`u $m_i$ est dans $L_i^\perp=L_i$.
Le vecteur de Wu $\omega_1-m_1=\omega_2+m_2$
convient.

Soit $q$ un repr\'esentant anisotrope du noyau de ${\rm Parf\/}$. D'apr\`es
le ${\petcap Lemme 1\/}$ un $b\/$-Lagrangien $L$ a un vecteur
de Wu $\omega$ avec $q\,(\omega)=0$ donc $\omega=0$ puisque
$q$ est anisotrope. Ainsi $q\,(L)=0$
donc $L=0$ et $V=0$, d'o\`u l'injectivit\'e.

Pour tout \'el\'ement
$\lambda$ de $K$,
le $b\/$-Lagrangien $K\,e$ de la forme quadratique $q_\lambda$
d\'efinie sur l'espace
vectoriel de base $e,\,f$ par $q_\lambda(e)=1,\, b_\lambda\,(e,f)=1$ et
$q_\lambda\,(f)=\lambda$
a $f$ comme vecteur de Wu et donc ${\rm Parf\/}\,(q)=[\lambda]$, d'o\`u
la surjectivit\'e.  
{\parindent=0pt\par\vskip .3cm
\vskip 0mm plus -20mm minus 1,5mm\penalty-50
{\bf 4.\ \/}%
{\petcap L'invariant de Arf\/}{\bf \pointir}%
\nobreak\parindent=20pt}%
Soit $k$ un corps de caract\'eristique $2$ dont une cl\^oture parfaite
est not\'ee $\iota : k\subset K$.
D'apr\`es la {\petcap Proposition\/} et le {\petcap Lemme 0\/} il y a
un morphisme de groupe ${\rm Arf\/} :WQ\,(k)\rightarrow k/{\cal P\/}\,(k)$
rendant commutatif le diagramme :
\[
  \begin{array}{ccccc}
    WQ\,(k) && \xrightarrow[]{\rm Arf\/} && k/{\cal P\/}\,(k) \\
    \ \ \downarrow{WQ\,(\iota)} &&&& \downarrow{\iota} \\
   WQ\,(K) && \xrightarrow[\sim]{\rm Parf\/ } && K/{\cal P\/}\,(K)\\
  \end{array}
\]
%$$\diagram{&\hdfl{\rm Arf\/}{}&k/{\cal P\/}\,(k)\cr
%\vbfl{\displaystyle }{}&&\vbfl{\displaystyle [\iota]}
%{\displaystyle \wr\kern-2pt}\cr
%			&\hdfl{\rm Parf\/}{\sim}&K/{\cal P\/}\,(K)\cr
%}$$
Soit $q : V\rightarrow k$ 
une forme quadratique non d\'eg\'en\'er\'ee sur $k$.
Sa forme bilin\'eaire $b$ a une base symplectique
$e_1, f_1,\ldots,e_n,f_n$. Alors, comme
$b\,(e_i,\, e_j)=b\,(f_i,\,f_j)=0$ et $b\,(e_i,\, f_j)=\delta_{i\,j}$
pour $i,\, j=1,\ldots,n$, 
l'\'el\'ement $\omega=\sqrt{q\,(e_1)}\,f_1+\cdots+\sqrt{q\,(e_n)}\,f_n$ 
de $V\otimes K$ est vecteur de Wu du Lagrangien
$K\,e_1\oplus\cdots\oplus K\,e_n$ de l'extension de $q$ \`a $K$ et
$q\,(\omega)=q\,(e_1)\,q\,(f_1)+\cdots +q\,(e_n)\,q\,(f_n)$. On retrouve
et compl\`ete le
\Thc Th\'eor\`eme de Arf
     {\rm (C.f. {\bf [1]} Satz 5, {\bf [3]} Appendix 1)}|
     La classe modulo ${\cal P\/}\,(k)$ de
$q\,(e_1)\,q\,(f_1)+\cdots +q\,(e_n)\,q\,(f_n)$
ne d\'epend pas du choix de la base symplectique.
C'est l'{\it invariant de Arf\/} ${\rm Arf\/}\,(q)$ de la forme quadratique
$q$, un morphisme
du groupe de Witt quadratique de $k$ sur le conoyau de son Artin-Schreier.
$${\rm Arf\/} : WQ\,(k)\rightarrow k/{\cal P\/}\,(k)
\xleftarrow[\sim]WQ\,(K)$$
s'identifiant au passage au groupe de Witt quadratique de la cl\^oture 
parfaite $K$.
\finc

{\parindent=0pt\par\vskip .3cm
\vskip 0mm plus -20mm minus 1,5mm\penalty-50
{\it Remerciements \/}{\bf \pointir}
Manifestons notre sympathie \`a
Monsieur Carreira dont les doigts habiles sont responsables de l'avenante
frappe du manuscript.}
\vskip .5 cm

\centerline{\bf R\'ef\'erences}
{\hangindent=1cm\hangafter=1\noindent{\bf [1]}\
{\bf C. Arf}\pointir{\sl Untersuchungen \H{u}ber quadratische Formen
in K\H{o}rpern der Characteristic 2.},
J. Reine Angew. Math. {\bf 183}, 148-167 (1941).

{\hangindent=1cm\hangafter=1\noindent{\bf [2]}\
{\bf J. Barge, J. Lannes, F. Latour {\rm et\/} P. Vogel}\pointir {\petcap Appendice}
de {\sl $\Lambda\/$-Sph\`eres},
Ann. scient. \'Ec. Norm. Sup., {\bf 8}, fasc. 4, 494-505 (1974).

{\hangindent=1cm\hangafter=1\noindent{\bf [3]}\
{\bf J. Milnor {\rm et\/} D. Husemoller}\pointir 
{\sl Symmetric bilinear forms},
Springer Verlag, (1973).

\vfill
\noindent{
Metteur en sc\`ene et secr\'etaire~: Alexis Marin Bozonat\hfill\break
\null\ courriel~: alexis.charles.marin@gmail.com}

%\vskip5mm
\noindent{
Institut Fourier, UMR 5582, Laboratoire de Math\'ematiques Universit\'e  Grenoble Alpes, CS 40700, 38058 Grenoble cedex 9, France}

\end{document}